\documentclass{amsart}
\usepackage{amssymb, amscd, enumerate}

\usepackage[matrix,arrow,curve,frame]{xy}
\xymatrixcolsep{1.9pc}
\xymatrixrowsep{1.9pc}
\newdir{ >}{{}*!/-5pt/\dir{>}}

\newcommand{\margnote}[1]{}

\newtheorem{main}{Theorem}

\newtheorem{theorem}{Theorem}[section]

\newtheorem{lemma}[theorem]{Lemma}
\newtheorem{corollary}[theorem]{Corollary}
\newtheorem{cor}[theorem]{Corollary}
\newtheorem{proposition}[theorem]{Proposition}
\newtheorem{prop}[theorem]{Proposition}
\newtheorem{definition}[theorem]{Definition}
\newtheorem{defn}[theorem]{Definition}

\newtheorem{numbered equation}[theorem]{(\hspace*{-0.1cm}}

\theoremstyle{definition}
\newtheorem{example}[theorem]{Example}
\newtheorem{remark}[theorem]{Remark}

\DeclareMathOperator{\Hom}{Hom}
\DeclareMathOperator{\HomC}{Hom_{\C}}

\DeclareMathOperator{\CHom}{Hom_{R}}

\DeclareMathOperator{\RCHom}{RHom_R}

\DeclareMathOperator{\mEnd}{{\underline{End}}}
\DeclareMathOperator{\End}{End}
\DeclareMathOperator{\EndR}{End_R}
\DeclareMathOperator{\EndC}{End_{\C}}

\DeclareMathOperator{\Mod}{Mod-}
\DeclareMathOperator{\Proj}{Proj-}
\DeclareMathOperator{\proj}{proj-}
\DeclareMathOperator{\modh}{mod-}

\DeclareMathOperator{\Ho}{Ho}

\DeclareMathOperator{\Ch}{Ch}

\newcommand{\cWE}{{\mathcal W}}

\newcommand{\bT}{{\mathbb T}}

\newcommand{\bZ}{{\mathbb Z}}
\newcommand{\uA}{\underline{A}}
\newcommand{\uR}{\underline{R}}
\newcommand{\uS}{\underline{S}}

\newcommand{\cA}{{\mathcal A}}
\newcommand{\cB}{{\mathcal B}}
\newcommand{\cC}{{\mathcal C}}

\newcommand{\DR}{{\mathcal D}_R}
\newcommand{\DS}{{\mathcal D}_S}
\newcommand{\DE}{{\mathcal D}_{\End_R(P)}}

\newcommand{\cE}{{\mathcal E}}
\newcommand{\cL}{{\mathcal L}}
\newcommand{\cM}{{\mathcal M}}
\newcommand{\cMc}{{\mathcal M}_c}
\newcommand{\cNc}{{\mathcal N}_c}
\newcommand{\cN}{{\mathcal N}}
\newcommand{\cP}{{\mathcal P}}

\newcommand{\cR}{{\mathcal R}}
\newcommand{\cS}{{\mathcal S}}
\newcommand{\cT}{{\mathcal T}}
\newcommand{\cU}{{\mathcal U}}
\newcommand{\cV}{{\mathcal V}}
\newcommand{\cW}{{\mathcal W}}

\newcommand{\cbU}{\overline{\cU}}

\newcommand{\C}{{\mathcal C}}
\newcommand{\D}{{\mathcal D}}
\newcommand{\Dc}{{\mathcal D}_c}

\newcommand{\U}{\mathcal U}

\newcommand{\EP}{\mathcal E(\mathcal P)}
\newcommand{\ab}{\mathcal Ab}
\newcommand{\EPo}{\End_R(P)}

\newcommand{\iso}{\cong}
\newcommand{\triequiv}{\he_{\Delta}}

\newcommand{\tens}{\otimes}

\newcommand{\mc}{\colon \,}

\newcommand{\ch}{{\mathcal C}h}
\newcommand{\chZ}{{\mathcal C}h_{\bZ}}

\newcommand{\ModA}{{\mbox {Mod-}}A}

\renewcommand{\to}{\longrightarrow}
\newcommand{\varrow}[1]{\hbox to #1{\rightarrowfill}}
\newcommand{\varl}[2]{\stackrel{#2}{\hbox to #1{\leftarrowfill}}}

\newcommand{\varrx}[2]{\stackrel{#2}{\hbox to #1{\rightarrowfill}}}
\newcommand{\parallelarrows}[1]{\begin{array}{c} {\hbox to
#1{\rightarrowfill}}  \vspace{-0.35cm} \\ {\hbox to
#1{\rightarrowfill}} \end{array}}

\newcommand{\he}{\simeq}
\newcommand{\dfn}{\textbf}
\newcommand{\adjoint}{\rightleftarrows}
\newcommand{\ovcat}{\!\downarrow\!}
\newcommand{\sSet}{s{\mathcal Set}}

\def\longfib{\DOTSB\relbar\joinrel\twoheadrightarrow}
\newcommand{\cofib}{\rightarrowtail}
\newcommand{\fib}{\twoheadrightarrow}
\newcommand{\trfib}{\stackrel{{\sim}}{\longfib}}
\newcommand{\trcofib}{\stackrel{{\sim}}{\rightarrowtail}}
\newcommand{\we}{\stackrel{{\sim}}{\longrightarrow}}

\newcommand{\ra}{\rightarrow}
\newcommand{\la}{\leftarrow}
\newcommand{\inc}{\hookrightarrow}
\newcommand{\wdo}{{\bullet}}  % The `dot' in S-dot construction
\newcommand{\Kpb}[1]{K_{+,hb}(\proj #1)}

\begin{document}

\title[Derived equivalences]{K-theory and derived equivalences} 

\date{September 6, 2002; 1991 AMS Math.\ Subj.\ Class.: 19D99, 18E30, 55U35}

\author{Daniel Dugger}
%\thanks{Department of Mathematics,  University of Oregon, 
%Eugene, OR 97403, USA
%(ddugger@math.uoregon.edu)}
\address{Department of Mathematics \\ University of Oregon \\
Eugene, OR 97403\\ USA}
\email{ddugger@math.uoregon.edu}

\author{Brooke Shipley}
%\thanks{Department of Mathematics,  Purdue University,  
%W. Lafayette, IN 47907,  USA
%(bshipley@math.purdue.edu)}
\address{Department of Mathematics \\ Purdue University\\  
W. Lafayette, IN 47907 
\\ USA}
\email{bshipley@math.purdue.edu}
\thanks{Second author supported in part by an NSF grant.}

\begin{abstract}
We show that if two rings have equivalent derived
categories then they have the same algebraic $K$-theory.
Similar results are given for $G$-theory, and for a large class of
abelian categories.
\end{abstract}

\maketitle

\tableofcontents

\section{Introduction}\label{sec-intro}

Algebraic $K$-theory began as a collection of elaborate invariants for
a ring $R$.  Quillen \cite{quillen} constructed these by feeding the
category of finitely-generated projective $R$-modules into  
the so-called Q-construction.  In fact, the Q-construction can
take as input any category with a sensible notion of exact
sequence.  Waldhausen later realized in~\cite{waldhausen} that the
same kind of invariants can be defined for a very broad class of {\it
homotopical\/} situations (Waldhausen used `categories with
cofibrations and weak equivalences').
To define the algebraic $K$-theory of a ring using the
Waldhausen approach, one takes as input the category of {\it bounded
chain complexes\/} of finitely-generated projective modules.

As soon as one understands this perspective it becomes natural to ask
whether the Waldhausen $K$-theory construction really depends on the
whole input category or just on the associated {\it homotopy
category\/} (where the weak equivalences have been inverted).  Or, in
the algebraic case, one asks whether the $K$-theory of a ring depends
only on the the associated derived category.  In this paper we answer
the latter question in the affirmative; if one is given the derived
category of a ring, together with its triangulation---but without
knowing which ring it is---then it is theoretically possible to
recover the algebraic $K$-theory of the ring.

\smallskip

We now give a more detailed description of the results.  If $R$ is a
ring, let $\DR$ denote the derived category of unbounded chain
complexes of $R$-modules.  Recall that $\DR$ is a triangulated
category in a standard way~\cite[10.4]{weibel}.  Also, let $K_*(R)$
denote the algebraic $K$-groups of $R$.  Our first theorem is the
following:

\begin{main}
\label{thm-K}
If $R$ and $S$ are two rings for which $\DR$ and $\DS$ are
equivalent as triangulated categories, 
then their algebraic $K$-groups are isomorphic: $K_*(R)\cong K_*(S)$.  
\end{main}

\noindent
When the hypothesis of the theorem holds we say that $R$ and $S$ are
\dfn{derived equivalent}, and so the result says that derived
equivalent rings have isomorphic $K$-theories.  (This definition of
`derived equivalent' is not manifestly the same as that of \cite[Def
6.5]{rickard-Morita}, but they do in fact
agree---see Theorem~\ref{thm-tilting one}).

\medskip

We can actually state somewhat stronger results.  Recall that $K_n(R)$
is the $n$th homotopy group of a certain space $K(R)$ produced by
ones favorite $K$-theory machine.  Let $\Dc(R)$ denote the full
subcategory of $\DR$ consisting of the perfect complexes---that is,
those complexes which are isomorphic in $\DR$ to a bounded complex of
finitely-generated projectives.  (The `c' is for `compact', a term
which is defined in Example~\ref{ex:compact}).

\begin{main}
\label{thm-Kequiv}
If $R$ and $S$ are rings such that $\Dc(R)$ and $\Dc(S)$ are
equivalent as triangulated categories, then $K_n(R)\cong K_n(S)$ for
all $n\geq 0$.  Even more, one has a weak equivalence of $K$-theory
spaces $K(R)\he K(S)$.  
\end{main}

The $K_0$ part of this result is very simple (see
\cite[9.3]{rickard-Morita}), and so our contribution is the extension
to higher $K$-theory.  We should mention that one can even weaken the
hypotheses somewhat, to require only an equivalence between $\Dc(R)$
and $\Dc(S)$ which commutes with the shift or suspension functor; see
Remark~\ref{rem-less}.

There are similar results for the $G$-theory of a ring.  Recall that
when $R$ is Noetherian $G(R)$ is the Quillen $K$-theory of the
category of finitely-generated $R$-modules (as opposed to
finitely-generated projectives); see \cite[Chapter 5]{Sr}.  In terms
of the Waldhausen machinery, it is the algebraic $K$-theory of the
category of bounded chain complexes of finitely-generated
$R$-modules---we denote the associated homotopy category by $\D_b(\modh
R)$.

\begin{main} Suppose that $R$ and $S$ are Noetherian rings.
\label{thm-Gequiv}
\begin{enumerate}[(a)]
\item
If $R$ and $S$ are derived equivalent, then $G(R)\he G(S)$; in
particular, $G_n(R)\cong G_n(S)$ for all $n\geq 0$.  
\item If $\D_b(\modh R)$ is triangulated-equivalent to $\D_b(\modh
S)$, then $G(R)\he G(S)$.
\end{enumerate}
\end{main}

Results along these lines first appeared in the work of Neeman
\cite{neeman}.  Neeman has had the much more ambitious goal of
actually constructing the algebraic $K$-theory space directly from the
derived category.  It seems he has accomplished this in the case of
abelian categories (cf. \cite[Thm. 7.1, p336]{neeman}), and so for
instance can construct $G(R)$ from $\D_b(\modh R)$ when $R$ is
Noetherian.  Using this result Neeman is able to prove
Theorem~\ref{thm-Gequiv}(b), and from this he is able to deduce
Theorem~\ref{thm-Kequiv} in the case of regular rings (because for
regular rings one has $G_*(R)\iso K_*(R)$).  Theorem~\ref{thm-Kequiv}
in the above generality is new, however, as are the other results
above.

Neeman's work is quite long and intricate, and it has sometimes been
met with a certain amount of suspicion---mostly because experts just
did not believe that $K$-theory could depend only on the derived
category.  The point we would like to accentuate is that our proofs of
the above theorems are all quite simple. 
The only `new' tool which enters the mix is the use of model categories.
Although model categories are not often used in these contexts, their
use effectively streamlines our work.  There are two main points
underlying the above theorems:

\medskip

{\em 
\begin{enumerate}[(1)]

\item Any equivalence of model categories yields a weak equivalence of
$K$-theory spaces (see Proposition~\ref{pr:QuilKiso}), and\\

\item If two rings $R$ and $S$ are derived equivalent then tilting
theory shows that their model categories of chain complexes $\Ch_R$
and $\Ch_S$ are in fact equivalent as model categories (see
Theorem~\ref{thm-tilting one}).
\end{enumerate}}

\medskip

\noindent
The first observation can be seen as an improvement
of~\cite[1.9.8]{thomason}, see Remark~\ref{rem-thom}.
The second is a more structured version of~\cite[6.4]{rickard-Morita}
and~\cite[3.3, 5.1]{rickard2}; note that unlike \cite{rickard2}, we do
not require any flatness hypotheses.  See also~\cite[5.1.1, B.1]{ss2}.

The observation in (2) is definitely surprising, although it turns out
that it is not hard to prove (in fact, considering the extra structure
in the model category seems to simplify the classical tilting theory
proofs).  The reason it is surprising is that the derived category of
$R$ is the `homotopy category' of $\Ch_R$, and this usually represents
only first-order information in the model category.  Equivalent model
categories have equivalent homotopy categories, but it almost never
works the other way around.  So something special happens when dealing
with chain complexes over a ring; the first order information here
determines all of the higher order information. Note that this does
not happen in arbitrary `abelian' model categories.  See also
Remarks~\ref{rem-Qeq} and~\ref{rem-marco}.

\smallskip

We state one last theorem along these lines, where we replace
the category of $R$-modules by any rich enough abelian category.  Of
course any abelian category $\cA$ has an unbounded derived category
$\D_\cA$, and we'll say that $\cA$ and $\cB$ are derived equivalent if
$\D_\cA$ is triangulated-equivalent to $\D_\cB$.  Let $K_c(\cA)$
denote the Waldhausen $K$-theory of the compact objects in $\ch_\cA$.
It turns out that the space $K_c(\Mod R)$ is just $K(R)$.

Recall that if $\cA$ is an abelian category, we say that an object $P$
is a {\em strong generator} if $X=0$ whenever $\hom_{\cA}(P,X) = 0$;
when $\cA$ has arbitrary coproducts, the object $P$ is called {\em
small} if $\oplus_{\alpha}\hom_{\cA}(P, X_{\alpha}) \to \hom_{\cA}(P,
\oplus_{\alpha} X_{\alpha})$ is a bijection for every set of objects
$\{X_{\alpha}\}_\alpha$.  Gabriel~\cite[V.1]{gabriel} has classified
the abelian categories which are equivalent to categories of modules
over a ring: these are the co-complete abelian categories with a
single strong generator. Freyd~\cite[5.3H]{freyd} generalized this to
include the case of many generators; see Theorems~\ref{thm-gab} and
\ref{thm-freyd}.  Using these basic tools, we can extend our above
statements to prove the following:

\begin{main}\label{thm-abresult}
Let $\cA$ and $\cB$ be co-complete abelian categories which have sets
of small, projective, strong generators.  Then
\begin{enumerate}[(a)]
\item $\cA$ and $\cB$ are derived equivalent if and only if 
$\ch_\cA$ and $\ch_\cB$ are equivalent as model categories.
\item If $\cA$ and $\cB$ are
derived equivalent, then $K_c(\cA) \he K_c(\cB)$.  
\end{enumerate}
\end{main}

Neeman \cite[7.1]{neeman} has proven that if $\cA$ and $\cB$ are small
abelian categories for which $\D_b(\cA)$ is triangulated-equivalent to
$\D_b(\cB)$, then $K(\cA) \he K(\cB)$ where $K(\cA)$ denotes the
Quillen $K$-theory of the exact category $\cA$.  There is little
overlap between this result and the above one: the abelian categories
in Theorem~\ref{thm-abresult} have infinite direct sums, so it follows
from the Eilenberg-Swindle that $K(\cA)$ and $K(\cB)$ are both trivial.
We do not know how to apply our methods to the kinds of abelian
categories Neeman deals with.

One final note: The reader may have noticed that we have always talked
about $K$-theory {\it spaces\/}, rather than $K$-theory {\it
spectra\/}.  In fact, all of the results in this paper hold when
restated in terms of spectra, and there is no difference in the
proofs.  We have chosen to avoid the added complications in an attempt
to streamline the presentation.

\subsection{Organization}
The proofs of Theorems A--C are given in Section~\ref{sec-proofs}, and
the paper has been structured so that the reader can get to them as
soon as possible.  The sections previous to that build up the
necessary machinery, but with most of the technical proofs postponed
until later.  Section~\ref{se:Kthy} recasts Waldhausen $K$-theory as
an invariant for model categories, and proves that it is preserved by
Quillen equivalences. Section~\ref{se:tilt} explains what tilting
theory has to say about Quillen equivalences between model categories
of chain complexes.  Finally, in Section~\ref{sec-many} we develop
the many-generators version of tilting theory, and prove Theorem D.

\subsection{Notation and terminology}

Being topologists, our convention is to always work with {\it chain\/}
complexes $C_*$ rather than cochain complexes.  So the differentials
have the form $d\mc C_n \ra C_{n-1}$, and the shift operator is
denoted as $\Sigma C$: it is the chain complex with $(\Sigma
C)_n=C_{n-1}$.

Throughout this paper we deal with right modules (and our rings are
not necessarily commutative).  Everything could be translated to left
modules as well, but because of the usual conventions for composing
maps, right modules are what naturally arise in some of our results;
see Theorems~\ref{thm-gab} and~\ref{thm-gabriel one} for example.
$\Mod R$ denotes the category of all $R$-modules, whereas $\modh R$
denotes the category of finitely-generated $R$-modules (we only use
this when $R$ is right-Noetherian).  Likewise, $\Proj R$ is the
category of all projective $R$-modules and $\proj R$ is the
subcategory of finitely-generated projectives.

Finally, if $\cC$ is a category then we write $\cC(X,Y)$ for
$\Hom_\cC(X,Y)$.  

\subsection{Acknowledgments}
We are grateful to Mike Mandell and Stefan Schwede
for several helpful conversations related to this paper.

%%%%%%%%%%%%%%%%%%%%%%%%%%%%%%%%%%%%%%%%%%%%%%%%%%%%%%%%%%%%%%%%%

\section{Model category preliminaries}

A model category is a category equipped with certain extra structures
which allow one to `do homotopy theory'.  The theory is based on three
standard examples: the category of topological spaces, the category of
simplicial sets, and the category of chain complexes over a given
ring.  In this section we recall the basic axioms of model categories,
and state the main facts we need in the body of the paper. \cite{DS}, 
\cite{hirschhorn} and \cite{hovey-book} are good references for
this material.

\medskip

\begin{defn}{\em
A \dfn{model category} is a category $\cM$ equipped with three
distinguished classes of maps: the {\it weak equivalences\/}, the {\it
cofibrations}, and the {\it fibrations\/}.  Cofibrations are depicted
as $\cofib$, fibrations as $\fib$, and weak equivalences as $\we$.
Maps which are both cofibrations and weak equivalences are called {\it
trivial cofibrations}, and denoted by $\trcofib$; {\it trivial
fibrations\/} are defined similarly.  The following axioms are
required:
\begin{enumerate}[{Axiom} 1:]
\item $\cM$ is complete and co-complete.
\item (Two-out-of-three axiom) If $f\colon A\ra B$ and $g\colon B \ra
C$ are maps in $\cM$ and any two of $f$, $g$, and $gf$ are weak
equivalences, then so is the third.
\item (Retract axiom) A retract of a weak equivalence
(respectively cofibration, fibration) is again a weak equivalence
(respectively cofibration, fibration).
\item (Lifting axiom) Suppose
\[
\xymatrix{ A \ar[r]\ar@{ >->}[d] & X \ar@{->>}[d] \\
           B \ar[r] & Y
}
\]
is a square (in which $A\ra B$ is a cofibration and $X\ra Y$ is a
fibration).  Then if either of the two vertical maps is a weak
equivalence, there is a lifting $B\ra X$ making the diagram commute.
\item\label{f} (Factorization axiom)
Any map $A\ra X$ may be functorially factored in two ways, as $A
\trcofib B \fib X$ and as $A \cofib Y \trfib X$.  
\end{enumerate}}
\end{defn}

Suppose maps $A\ra B$ and $X\ra Y$ are given. When any square
as in Axiom 4 has a lifting $B\ra X$, we say that $A\ra B$ has the
{\it left-lifting-property\/} with respect to $X\ra Y$.

\begin{example}  In this paper we only deal explicitly with model categories 
on categories of chain complexes.
\begin{enumerate}[(a)]
\item
The category $\ch_R^+$ of non-negatively graded chain complexes over a
ring $R$ has a model structure where 
the weak equivalences are the maps inducing homology isomorphisms (the
\dfn{quasi-isomorphisms}), the fibrations are the maps which are
surjective in positive degrees, and the cofibrations are the
monomorphisms with degreewise projective cokernels;  
see ~\cite[II p.\ 4.11, Remark 5]{Q}, ~\cite[Sec. 7]{DS}.
This model
structure on $\ch_R^+$ is referred to as the {\it projective\/} model
structure since there are other model structures on $\ch_R^+$.

\item
The category $\ch_R$ of unbounded chain complexes over a ring $R$
also has a {\bf (projective) model structure} with weak equivalences
the homology isomorphisms, and fibrations the epimorphisms;
see~\cite[2.3.11]{hovey-book},~\cite[5]{ss1}.  Every cofibration is
still a degreewise split injection and the cokernel is levelwise
projective,
but not all such degreewise split injections are cofibrations.

\item $\ch_R$ has another model structure with the same weak
equivalences, but where the cofibrations are the monomorphisms.
The fibrations are harder to describe, but any fibration is a 
degreewise surjection with levelwise-injective kernel.  This is 
the {\em injective model structure} on $\ch_R$.  In this paper we
only need to use the projective model structure on $\ch_R$,
however.
\end{enumerate}
\end{example}

When $\cM$ is a model category, one may formally invert the weak
equivalences $\cWE$ to obtain the category-theoretic localization
$\cWE^{-1}\cM$.  This is the \dfn{homotopy category} of
$\cM$, written $\Ho \cM$; see~\cite[I.1]{Q},~\cite[6.2]{DS}.  Since
the weak equivalences in $\ch_R$ are the quasi-isomorphisms, the
homotopy category $\Ho \ch_R$ is equivalent to the (unbounded) derived
category $\DR$ (cf. \cite[Example 10.3.2]{weibel}).  

A model category is called \dfn{pointed} if the initial object and
terminal object are the same.  The homotopy category of any pointed
model category turns out to have a {\it suspension functor\/}
$\Sigma$.  For topological spaces this is ordinary suspension, whereas
for $\ch_R^+$ and $\ch_R$ it is the functor sending a chain complex
$C$ to the shift $\Sigma C$ with $(\Sigma C)_n = C_{n-1}$.  As the
example of $\ch_R^+$ shows, this functor need not be an equivalence.
When it {\it is\/} an equivalence we say that $\cM$ is a \dfn{stable}
model category, and in this case $\Ho \cM$ becomes a triangulated
category in a natural way~\cite[7.1]{hovey-book}.  (When $\cM$ is not
stable, $\Ho \cM$ only has a `partial' triangulation; see~\cite[I.2,
I.3]{Q},~\cite[6.5]{hovey-book} for details).  For $\ch_R$ this of
course specializes to the usual triangulation on $\DR$.

\begin{definition} \label{def-Quillen}
{\em A \dfn{Quillen map} of model categories $\cM \ra \cN$ consists of a
pair of adjoint functors $L\colon \cM \adjoint \cN \colon R$ such that
$L$ preserves cofibrations and trivial cofibrations (it is equivalent
to require that $R$ preserves fibrations and trivial fibrations).  
In this case the pair $(L,R)$ is also called a {\it Quillen pair\/}.}
\end{definition}

\begin{example}
\label{ex:Quilpair}
Let $R\ra S$ be a map of rings.  The adjoint pair of functors
$\cL\colon \Mod R \adjoint \Mod S \colon \cR$ defined by
$\cL(M)=M\tens_R S$ and $\cR(N)=\Hom_R(S,N)$ prolongs to an adjoint
pair between categories of chain complexes.  One readily checks that
these prolongations are 
Quillen maps $\ch_R^+ \ra \ch_S^+$ and $\ch_R \ra \ch_S$.
\end{example}

A Quillen map induces adjoint total derived functors between the
homotopy categories~\cite[I.4]{Q}.  The map is a \dfn{Quillen
equivalence} if the total derived functors are adjoint equivalences of
the homotopy categories.  This is equivalent to Quillen's original
definition by~\cite[1.3.13]{hovey-book}.  More generally we say that
$\cM$ and $\cN$ are \dfn{Quillen equivalent} if they are connected by a
zig-zag of Quillen equivalences, and we write $\cM\he_Q \cN$.
As one simple example, the identity functors give a Quillen equivalence
$\ch_R^{\text{proj}} \ra \ch_R^{\text{inj}}$ between the projective
and injective model structures on $\ch_R$.  

\begin{remark}\label{rem-Qeq}
In general, having a Quillen equivalence of model categories is much
stronger than just having an equivalence between the associated
homotopy categories.  This is because of the added structure required
for a Quillen map; functors on the homotopy categories may not lift to
the model category level, and even if they do they may not be
compatible with the model category structures.  For example, it
follows from~\cite[I.4 Thm.~3]{Q} that Quillen maps between stable
model categories induce triangulated functors between the homotopy
categories.  Quillen maps preserve even more structure, for example
the simplicial mapping space
structures~\cite[5.4]{DK},~\cite[5.6.2]{hovey-book}.  There are simple
topological examples---see \cite[3.2.1]{ss2}, for instance---of stable
model categories which have the same triangulated homotopy category,
but which are nevertheless not Quillen equivalent.  In
Remark~\ref{rem-marco} we discuss another example (based
on~\cite{marco}) which is entirely algebraic.
\end{remark}

The following theorem shows that for the special case of the model
categories $\ch_R$, Quillen equivalence is not a stronger notion than
triangulated equivalence of homotopy categories.  In some sense this
happens because rings are determined by `first order'
information---compared, for example, to differential graded rings
which are not.  This is proved in Section~\ref{sec-proof de iff qe} as
Theorem~\ref{thm-tilting one} (Parts ~\ref{1} and~\ref{2}).

\begin{theorem}\label{thm-de iff qe}
Two rings $R$ and $S$ are derived equivalent if and only if
their associated model categories of chain complexes $\ch_R$ and $\ch_S$
are Quillen equivalent.
\end{theorem}

This theorem cannot be extended to cover the case where $R$ or $S$ is
a differential graded algebra; we give an example in \cite{ds-marco}
which is discussed a little in Remark~\ref{rem-marco}.  In
Corollary~\ref{cor-abcat} we do give a certain extension of this
theorem to abelian categories, however.  The situation is a little
confusing, because these two sentences may seem
contradictory. They are not though; see Remark~ \ref{caution}.

%%%%%%%%%%%%%%%%%%%%%%%%%%%%%%%%%%%%%%%%%%%%%%%%%%%%%%%%%%%%%%%%%%%%%%%%%%

\section{$K$-theory and model categories}
\label{se:Kthy}

In \cite[Section 1]{waldhausen} Waldhausen defined a notion of {\it
category with cofibrations and weak equivalences} and showed how to
construct a $K$-theory space from such data.  The purpose of this
section is to adapt Waldhausen's machinery to the context of model
categories.  This is almost entirely straightforward, but it has the
advantage of streamlining the theory somewhat.

\medskip

Let $\cM$ be a pointed model category with initial object $*$.  An
object $A$ is called \dfn {cofibrant} if $* \cofib A$ is a
cofibration.  By a \dfn{Waldhausen subcategory} of $\cM$ we mean
a full subcategory $\cU$ with the properties that
\begin{enumerate}[(i)]
\item $\cU$ contains the initial object $*$;
\item Every object of $\cU$ is cofibrant;
\item If $A\cofib B$ and $A\ra X$ are maps in $\cU$, then the pushout
$B\amalg_A X$ (computed in $\cM$) belongs to $\cU$.
\end{enumerate}

The proof of the following is just a matter of chasing through the
definitions:

\begin{lemma}
Any Waldhausen subcategory of $\cM$, equipped with
the notions of cofibrations and weak equivalences from $\cM$, is a
{\rm `category with cofibrations and weak equivalences'} in the sense
of \cite[1.2]{waldhausen}; also, it satisfies the saturation axiom
\cite[p. 327]{waldhausen}.
\end{lemma}

The lemma says that we may apply Waldhausen's $S_{\wdo}$-construction
\cite[1.3]{waldhausen} to obtain a simplicial category
$wS_\wdo(\cU)$.  Taking the nerve in every dimension gives a
simplicial space $[n]\mapsto N(wS_n(\cU))$, and $K(\cU)$ is defined to
be loops on the realization of this simplicial space: $K(\cU)=\Omega
|NwS_\wdo(\cU)|$.  One defines the algebraic $K$-groups of $\cU$ by
$K_n(\cU)=\pi_n(K(\cU))$.

We give a partial description of $wS_\wdo(\cU)$ here, because we 
need it later in Appendix A.  
Let $wF_n(\cU)$ denote the category whose objects are
sequences $\{A\}$ of cofibrations $A_0 \cofib A_1 \cofib \cdots
\cofib A_n$ in $\cU$, and whose morphisms are commutative diagrams
$\{A\} \ra \{A'\}$ in which every map $A_n \ra A'_n$ is a weak
equivalence.  One can almost make $[n]\mapsto wF_{n-1}(\cU)$ into a
simplicial category (where $wF_{-1}(\cU)$ is interpreted as the trivial
category with one object and an identity map) by defining
\[ d_i \bigl ( [A_0
\cofib A_1 \cofib \cdots \cofib A_n] \bigr) =
\begin{cases}
[A_0 \cofib \cdots \cofib \hat{A}_i \cofib \cdots \cofib A_n] &
 \text{\ if\ } i\neq 0,  \\
{[A_1 / A_0 \cofib A_2 / A_0 \cofib \cdots \cofib A_n / A_0]} &
 \text{\ if\ } i=0.
\end{cases}
\]
The difficulty is that with this definition the simplicial identities
do not hold on the nose, in the end because there are different
possible choices for the quotients $A_i/A_0$ (they are canonically
isomorphic, but still different).  The category $wS_n(\cU)$ is
equivalent to $wF_{n-1}(\cU)$, but is slightly `fatter' in a way that
allows one to make the face and degeneracy maps commute on the nose.
The reader is referred to \cite[p. 328]{waldhausen} for the precise
definition---it should be noted, though, that the basic ideas in the
present paper can all be understood by pretending that $wS_n(\cU)$ is
just $wF_{n-1}(\cU)$.  The only time the details of $wS_n(\cU)$ are
needed is in the Appendix.

\begin{example}
\label{ex:U for rings}
Let $R$ be a ring.
The following are Waldhausen subcategories of $\ch_R$ (as is easily
verified).
\begin{enumerate}[(1)]
\item $\cU_K=\{$all bounded complexes of finitely-generated projectives$\}$.
\item $\cU_G=\{$all bounded below complexes $C$ of finitely-generated
projectives such 
%DD: Brute force spacing.  Surely not the best way to do this...
\qquad\phantom{A}\qquad\phantom{A}
that $H_k(C)\neq 0$ for only
finitely-many values of $k\}$.
\end{enumerate}
\end{example}

Let $K(R)$ and $G(R)$ denote the Quillen $K$-theory spaces for the
exact categories of finitely-generated projectives and
finitely-generated modules, respectively.  Then we have:

\begin{lemma}
$K(\cU_K) \he K(R)$, and if $R$ is Noetherian then $K(\cU_G)\he G(R)$.
\end{lemma}

\begin{proof}
A reference for $K(\cU_K)\he K(R)$ is \cite[1.11.7]{thomason}.  For
the $G$-theory, the reference is \cite[3.11.10, 3.12,
3.13]{thomason}---however, since the terminology of that paper is
fairly cumbersome, we repeat the proof for the reader's
convenience.

Let $\cV$ denote the subcategory of $\ch_R$ consisting of all bounded
complexes of finitely-generated modules; \cite[1.11.7]{thomason} shows
that $K(\cV)$ is the same as $G(R)$.  Let $\cW$ denote the subcategory
of $\ch_R$ consisting of all chain complexes quasi-isomorphic to an
element of $\cV$.  One can check that if $R$ is Noetherian then
$\cU_G$ consists precisely of the cofibrant objects in $\cW$.  Then
\cite[1.9.8]{thomason} shows that $\cU_G \inc \cW$ and $\cV\inc \cW$
induce equivalences of $K$-theory spaces.  
\end{proof}

\begin{example}
\label{ex:compact}
If $\cT$ is a triangulated category with infinite sums, an object
$X\in \cT$ is called \dfn{compact} if the natural map
$\oplus_\alpha \cT(X,Z_\alpha) \ra \cT(X,\oplus_\alpha Z_\alpha)$ is a
bijection for every collection $\{Z_\alpha\in \cT\}$.  If $\cM$ is a
stable model category, it is easy to check that the homotopy category
$\Ho\cM$ has all infinite sums.  We'll say that an object in $\cM$ is
compact if its image in $\Ho\cM$ is compact.  The subcategory
$\cMc\subseteq \cM$ consisting of all compact, cofibrant objects is a
complete Waldhausen category.

We are especially interested in this for the case $\cM=\ch_R$, where a
theorem of B{\"o}kstedt-Neeman \cite[6.4]{boek-nee} identifies the
compact objects as the perfect complexes, i.e. the complexes which are
quasi-isomorphic to a bounded complex of finitely-generated
projectives.
\end{example}

\begin{example}
Waldhausen never explicitly used model categories, but he could have
been working in this context all along.  Waldhausen developed his
machinery to apply to the following case.  Let $X$ be a simplicial
set, and let $(X\ovcat \sSet \ovcat X)$ denote the category of
retractive spaces over $X$.  This has a natural model structure
inherited from the category of simplicial sets~\cite[II.3]{Q} by
forgetting the retraction over $X$ (cf. \cite[7.6.5]{hirschhorn}).
Take $\cU$ to be the subcategory consisting of those retractive spaces
$X\inc Z \ra X$ for which the map $X\inc Z$ is obtained by attaching
finitely many simplices.  This is a Waldhausen subcategory, and the
associated $K$-theory space is denoted $A(X)$;
see~\cite[2.1]{waldhausen}.
\end{example}

If $\cU$ is a subcategory of $\cM$, write $\cbU$ for the full
subcategory of $\cM$ consisting of all cofibrant objects which are
weakly equivalent to an object in $\cU$.  From
Example~\ref{ex:compact} above it follows that
$(\ch_R)_c=\overline{\cU_K}$ (where $\cU_K$ is from Example~\ref{ex:U
for rings}).  Call the Waldhausen category $\cU$ \dfn{complete} if
$\cU=\cbU$.

Suppose that $(L,R)\colon \cM \ra \cN$ is a Quillen map of pointed model
categories.  Let $\cU$ and $\cV$ be Waldhausen subcategories of $\cM$
and $\cN$ such that $L$ maps $\cU$ into $\cV$.  Since $L$ preserves
cofibrations, one checks easily that it induces a well-defined map
$K(\cU) \ra K(\cV)$.

\begin{prop}
\label{pr:QuilKiso}
Suppose that $(L,R)$ is a Quillen equivalence, and that $\cU$ is a
complete Waldhausen subcategory of $\cM$.  Let
$\cV=\overline{L\cU}$---i.e., $\cV$ consists of all cofibrant objects
which are weakly equivalent to an object in $L(\cU)$.  Then $\cV$ is a
complete Waldhausen subcategory of $\cN$, and $L\colon K(\cU) \ra
K(\cV)$ is a weak equivalence.
\end{prop}

The proof is simple but long winded, so we defer it to an appendix.

\begin{remark}
\label{rem-QuilKisoR}
The proposition also works in the following way.  Let $Q$ be a
cofibrant replacement functor for $\cM$; for example one can take the
map $* \ra X$ and apply the functorial factorization $* \cofib QX
\trfib X$ in $\cM$ to define $Q$.  Similarly, let $F$ be a
fibrant-replacement functor for $\cN$ with $Y \trcofib FY \fib *$ for
$Y$ in $\cN$.  Suppose that $\cV$ is a complete Waldhausen subcategory
of $\cN$.  Define $R\cV$ to be the set of all objects of the form
$QRFX$ where $X\in \cV$, and let $\cU=\overline{R\cV}$.  Then $\cU$ is
a complete Waldhausen subcategory of $\cM$, and $\overline{L\cU}=\cV$.
The functor $L$ induces a map $K(\cU) \ra K(\cV)$, and the proposition
says this is an equivalence.  So we have actually proven:
\end{remark}

\begin{cor}
\label{co:QuilKiso2}
Let $\cM$ and $\cN$ be model categories connected by a zig-zag of
Quillen equivalences.  Let $\cU$ be a complete Waldhausen subcategory
of $\cM$, and let $\cV$ consist of all cofibrant objects in $\cN$
which are carried into $\cU$ by the composite of the derived functors 
of the Quillen equivalences.  Then $\cV$ is a complete Waldhausen subcategory 
of $\cN$, and there is an induced zig-zag of weak equivalences between
$K(\cU)$ and $K(\cV)$.
\end{cor}

\begin{cor}
\label{co:stable QuilKiso}
A Quillen equivalence $\cM \ra \cN$ between stable model categories
induces a weak equivalence of $K$-theory spaces $K(\cMc)\we
K(\cNc)$, where $\cMc$ and $\cNc$ denote the subcategories of
cofibrant, compact objects.
\end{cor}

\begin{proof}
Write the functors of the Quillen equivalence as $(L,R)$.  
The derived functors of $L$ and $R$ induce an equivalence between the
homotopy categories, and so in particular they take compact objects to
compact objects.  This clearly implies $\cNc \supseteq
\overline{L\cMc}$; it basically gives the opposite inclusion as well,
but we now explain this in more detail.

If $X$ is in $\cNc$, let $F X$ be a fibrant-replacement 
$X \trcofib FX \fib *$ 
in $\cN$ and let $Q(R F X)$ be a cofibrant-replacement 
$* \cofib Q(R F X) \trfib R F X$ 
of $R F X$ in $\cM$.  Because the
derived functors of $L$ and $R$ take compact objects to compact
objects, $QR F X$ must still be compact---i.e., $QR F X \in
\cMc$.  Yet $L Q R F X$ is weakly equivalent to $X$, and so $X\in
\overline{L\cMc}$.  At this point we have shown $\cNc=
\overline{L\cMc}$, and so we can just apply
Proposition~\ref{pr:QuilKiso}.
\end{proof}

\begin{cor}\label{cor-QuilKiso}
If $\ch_R$ and $\ch_S$ are Quillen equivalent (perhaps through a
zig-zag of Quillen equivalences), then $K(R)\he K(S)$.
\end{cor}

\begin{proof}
We have already remarked that $K(R)\he K(\cU_K)$, and
$\cU_K=(\ch_R)_c$.  All the intermediate model categories in the
zig-zag must be stable because `stability' is preserved under Quillen
equivalence.  Therefore Corollary~\ref{co:stable
QuilKiso} applies.
\end{proof}

\begin{remark}\label{rem-thom}
The above corollary has two improvements over similar results in the
literature.  The first is that we are allowing a zig-zag of Quillen
equivalences, rather than just an equivalence $\ch_R\ra \ch_S$; in
particular, note that our zig-zag could conceivably pass through very
non-algebraic model categories.  For just a single Quillen equivalence
$\ch_R\ra \ch_S$, the closest result in the literature seems to be
\cite[1.9.8]{thomason}.  In that result, however, the functor
$L\colon \ch_R \ra \ch_S$ is required to be {\it complicial\/},
meaning in part that it is induced via prolongation from a functor
$\Mod R \ra \Mod S$.  In some of our applications $L$ is the
functor which tensors with a chain complex of projectives (rather than
just a single projective), and so the \cite{thomason} result is not
applicable.
\end{remark}

%%%%%%%%%%%%%%%%%%%%%%%%%%%%%%%%%%%%%%%%%%%%%%%%%%%%%%%%%%%%%%%%%%%%%%%%%%%%

\section{Tilting Theory}
\label{se:tilt}

In this section we determine 
the {\it algebraic content\/} of having a Quillen equivalence between
$\ch_R$ and $\ch_S$ for rings $R$ and $S$.  
A nice, complete answer can be given in terms of {\it tilting
theory\/}.  Originally tilting theory only dealt with derived
equivalences, but it turns out that for rings derived equivalence and
Quillen equivalence coincide.

\medskip

We begin with a classical analogue of tilting theory, namely Morita
theory.  Morita theory describes necessary and sufficient conditions
for when two categories of modules are equivalent.  Call a (right)
$R$-module $P$ a \dfn {strong generator} if $\hom_R(P,X)=0$ implies $X
=0$ for any (right) $R$-module $X$.

\begin{theorem}\label{thm-Morita}{\bf (Morita Theory)}
Given rings $R$ and $S$, the following conditions are equivalent:
\begin{enumerate}
\item\label{one} The categories of (right) modules over $R$ and $S$
are equivalent.
\item\label{two} There is an $R$-$S$ bimodule $M$ and an $S$-$R$
bimodule $N$ such that $M \otimes_S N \iso R$ as $R$-bimodules and $N
\otimes_R M \iso S$ as $S$-bimodules.
\item\label{three} There is a (right) $R$-module $P$ which is
finitely-generated, projective and a strong generator such that
$\hom_R (P,P) \iso S$.
\end{enumerate}
\end{theorem}

\begin{proof}
We only give a brief sketch because this is classical,
see~\cite[9.5]{weibel}.  For (\ref{two}) implies (\ref{one}), the
functors $- \otimes_R M \mc \Mod R \to \Mod S$ and $- \otimes_S N \mc
\Mod S \to \Mod R$ give the inverse equivalences.  For (\ref{one})
implies (\ref{three}), given an equivalence $F \mc \Mod S \to \Mod R$
one may take $P= F(S)$.  For (\ref{three}) implies (\ref{two}), take
$N = P$ since $P$ is a $\hom_R(P,P)$-$R$ bimodule and take $M =
\hom_R(P,R)$ which is an $R$-$\hom_R(P,P)$ bimodule.
\end{proof}

Now we turn to the analogue of Morita theory for categories of chain
complexes, called `tilting theory'.  This analogue was developed
by Rickard in~\cite[6.4]{rickard-Morita} to classify derived
equivalences of rings.  
Later, Keller~\cite[8.2]{dg-keller} broadened
tilting theory to apply to more general derived equivalences of
abelian categories.  We extend both sets of results to give Quillen
equivalences underlying the derived equivalences.
Theorem~\ref{thm-tilting one} below extends Rickard's work, whereas
the generalization to abelian categories is considered in
Section~\ref{sec-many}.  These results can also be used to remove
certain flatness assumptions in~\cite[3.3, 5.1]{rickard2}.

Let $\cT$ be a triangulated category.  Recall that a \dfn{full
subtriangulated category} $\cS$ is a full subcategory which is (i)
closed under isomorphisms, (ii) closed under the suspension functor,
and (iii) has the property that if two objects of a distinguished
triangle in $\cT$ lie in $\cS$ then so does the third object.  When
$\cT$ has infinite sums, a full subtriangulated category is called
\dfn{localizing} if it is closed under coproducts of sets of
objects~\cite[1.5.1, 3.2.6]{neeman-book}.  A complex $P$ in $\cT$ is a
(weak) \dfn{generator} if the only localizing subcategory of $\cT$
which contains $P$ is $\cT$.  Although this definition looks much
different than the definition of a strong generator, it is not.  If
$P$ is compact (see Example~\ref{ex:compact} for a definition), then
$P$ is a (weak) generator if and only if 
$\cT(P,X)_* = 0$ implies $X$ is trivial (see \cite[2.2.1]{ss2} for
a proof that these are equivalent).  Here $\cT(-,-)_*$
denotes the graded maps with $\cT(X,Y)_n = \cT(\Sigma^n X, Y)$.

An object $P\in \ch_R$ is called a \dfn{tilting complex} if it
is a bounded complex of finitely-generated
projectives, a generator of $\DR$, and $\DR(P,P)_*$ is
concentrated in degree zero \cite[Def. 6.5]{rickard-Morita}.  Here is
our generalization of Rickard's result \cite[Thm 6.4]{rickard-Morita}:

\begin{theorem} \label{thm-tilting one} {\bf (Tilting theorem)}
The following conditions are equivalent
for rings $R$ and $S$:   
\begin{enumerate}
\item\label{1}
There is a zig-zag of Quillen equivalences between 
the model categories of chain complexes of $R$- and $S$-modules:
\[\ch_R \he_Q \ch_S. \] 
\item\label{2}
The unbounded derived categories are triangulated equivalent:
\[\DR \triequiv \DS.\] 
\item\label{2'} 
The naive homotopy categories of bounded chain
complexes of finitely generated projective $R$ and $S$-modules are
triangulated equivalent:
\[ K_b(\proj R)  \triequiv
K_b(\proj S).\]
\item \label{3} The model category $\ch_R$ has a tilting complex 
$P$ whose endomorphism ring in $\DR$ is isomorphic to $S$:
$\DR(P,P)\iso S$. 
\end{enumerate}
\end{theorem}

\begin{remark}
Rickard \cite[6.4]{rickard-Morita} showed that $(3)$ and $(4)$ are
equivalent and that both these are equivalent to having a triangulated
equivalence $\D_b(\Mod R)\he \D_b(\Mod S)$.  He defined two rings to
be `derived equivalent' if any of these conditions holds.  We defined
`derived equivalent' to mean (2), and so the result shows that our use
agrees with Rickard's.  Note that \cite[6.4]{rickard-Morita}
gives several other equivalent conditions involving variations of the
derived category; see Proposition~\ref{pr-DRsubcat} as well.
\end{remark}

\begin{proof}[Proof of $(1)\Rightarrow (2) \Rightarrow (3) \Rightarrow (4)$]
Every Quillen equivalence of stable model categories induces an
equivalence of triangulated homotopy categories~\cite[I.4 Theorem 3]{Q}, 
so (\ref{1}) implies (\ref{2}).  Any triangulated equivalence restricts to
an equivalence between the respective subcategories of compact
objects.  Since $K_b(\mbox{proj-}R)$ is equivalent to the full
subcategory of compact objects in $\DR$ by \cite[Prop.\
6.4]{boek-nee}, (\ref{2}) implies (\ref{2'}).

Now we assume condition (\ref{2'}) and choose a triangulated
equivalence between $K_b(\mbox{proj-}R)$ and $K_b(\mbox{proj-}S)$.
Let $S[0]$ be the free $S$-module on one generator, viewed as a
complex in $\ch_S$ concentrated in dimension zero; let $T$ be its
image in $K_b(\proj R)$.  We have $\DR(T,T) \iso \DS(S[0],S[0]) \iso
S$.  Since $S[0]$ generates $K_b(\mbox{proj-}S)$, $T$ generates
$K_b(\mbox{proj-}R)$.  Since $R[0]$ is a generator of $\DR$ and
$R[0]\in K_b(\proj R)$, the only localizing subcategory of $\DR$
containing $K_b(\proj R)$ is $\DR$; so $T$ generates $\DR$.  Hence $T$
is a tilting complex and condition (\ref{3}) holds.
\end{proof}

The real content of the theorem, of course, is the proof that
$(4)\Rightarrow (1)$.  This is given in Section~\ref{sec-proof de iff
qe}, after we have developed a little more machinery.

\begin{remark}\label{rem-less}
We could have put one more intermediary condition in
Theorem~\ref{thm-tilting one}.  Instead of a triangulated equivalence
(in either (2) or (3)) we could have required only an equivalence of
categories which commutes with the shift or suspension functor.  Such
an equivalence would preserve compact objects and preserve the graded
maps $\D(-,-)_*$.  It would also preserve the property of being a
compact generator, since an object is a compact generator if and only
if it detects trivial objects by~\cite[2.2.1]{ss2}.  Thus, such
equivalences preserve tilting complexes.  We do not have very
interesting examples of such equivalences, though (other than
triangulated equivalences).
\end{remark}

\begin{remark}\label{rem-ss2}
The two tilting theory results in this paper, Theorem~\ref{thm-tilting
one} and its analogue Theorem~\ref{thm-tilting}, also appear in
disguised form in~\cite{ss2}.  Chain complexes do not satisfy the
stated hypotheses of the tilting theorem in~\cite[5.1.1]{ss2}, but
in~\cite[Appendix B.1]{ss2} chain complexes are shown to be Quillen
equivalent to a model category which does satisfy the stated
hypotheses.  So Theorems~\ref{thm-tilting one} and~\ref{thm-tilting}
can be considered as special cases of~\cite[5.1.1]{ss2}.
Here, though, we have removed all hypotheses and the proofs are much
simplified---they only use categories of chain complexes, whereas the
proofs in~\cite{ss2} require the use of the new symmetric monoidal
category of symmetric spectra~\cite{hss}.
\end{remark}

%%%%%%%%%%%%%%%%%%%%%%%%%%%%%%%%%%%%%%%%%%%%%%%%%%%%%%%%%%%%%%%%%%%%%%%%%%%%

%%%%%%%%%%%%%%%%%%%%%%%%%%%%%%%%%%%%%%%%%%%%%%%%%%%%%%%%%%%%%%%%%%%%%%%%%%%%

\section{Proofs of the main results}
\label{sec-proofs}

If you accept the basic results stated so far, it becomes easy
to prove the first three theorems cited in the introduction.  

\smallskip

\begin{proof}[Proof of Theorem~\ref{thm-Kequiv}]
This follows from Corollary~\ref{cor-QuilKiso} together with the
equivalence of Parts~\ref{1} and~\ref{2'} in Theorem~\ref{thm-tilting
one}.  Note that $\Dc(R)$ and $K_b(\proj R)$ are two names for the
same thing, by \cite[6.4]{boek-nee}.
\end{proof}

\begin{proof}[Proof of Theorem~\ref{thm-K}]
If $\DR$ and $\DS$ are equivalent as triangulated categories, then so
are their full subcategories of compact objects.  So
Theorem~\ref{thm-Kequiv} applies.  This also follows from
Corollary~\ref{cor-QuilKiso} and the equivalence of Parts~\ref{1}
and~\ref{2} in Theorem~\ref{thm-tilting one}.
\end{proof}

We now turn our attention to the proof of Theorem~\ref{thm-Gequiv},
which is the $G$-theory result.  We begin with a proposition which is
fairly interesting in its own right.  Consider a function $\cC$ which
assigns to each ring $R$ a subcategory of $\DR$.  We say that the
assignment \dfn{preserves equivalences} if every triangulated
equivalence $F\colon \DR \ra \DS$ restricts to an equivalence between
$\cC(R)$ and $\cC(S)$.

Here is some new notation: $\D_{h+}(\Mod R)$ denotes the full
subcategory of $\DR$ consisting of chain complexes with bounded below
homology, and $\D_{hb}(\Mod R)$ denotes the full subcategory of
complexes with bounded homology.  One can similarly define
$K_{hb}(\proj R)$, etc.  The notation $\Kpb{R}$ means the
intersection of $K_{+}(\proj R)$ and $K_{hb}(\proj R)$.  It is an easy
exercise to check that $\D_{h+}(\Mod R)=K_{+}(\Proj R)$ and
$\D_{hb}(\Mod R)=\D_b(\Mod R)$.

\begin{prop}
\label{pr-DRsubcat}
The assignments $R\mapsto \cC(R)$ preserve 
equivalences, where $\cC(R)$ is any of the following:
\[
K_b(\proj R), \quad
K_+(\Proj R)=\D_{h+}(\Mod R), \quad
\D_{h-}(\Mod R),
\]
\[
\D_{hb}(\Mod R)=\D_b(\Mod R), \quad
K_+(\proj R), \quad
\Kpb{R}.
\]
\end{prop}

\begin{proof}
The result \cite[6.4]{boek-nee} identifies $K_b(\proj R)$ with the
subcategory of compact objects in $\DR$.  Any equivalence $\DR \ra
\DS$ must preserve direct sums, and so it takes compact objects to
compact objects.

A complex $X$ lies in $\D_{h+}(\Mod R)$ if and only if it satisfies the
following property: for any compact object $A$, there exists an $N$
such that $\DR(\Sigma^{-k}A,X)=0$ for $k>N$.  Since triangulated 
equivalences preserve compact objects and the suspension, they 
preserve these objects as well.

Similarly, a complex $X$ lies in $\D_{h-}(\Mod R)$ if and only if for
any compact object $A$, there exists an $N$ such that $\DR(\Sigma^k
A,X)=0$ for all $k>N$.  The same argument as above applies.
For $\D_{hb}(\Mod R)$, note that this is just the intersection of
$\D_{h+}(\Mod R)$ and $\D_{h-}(\Mod R)$.

The case of $K_+(\proj R)$ is harder, but was proven by Rickard---see
the first paragraph in the proof of \cite[8.1]{rickard-Morita}.
Finally, $\Kpb{R}$ is just the intersection of $K_+(\proj R)$ and
$\D_{hb}(\Mod R)$.  
\end{proof}

$\Kpb{R}$ is the full subcategory of $\DR$ consisting of complexes
which are quasi-isomorphic to a bounded-below complex of
finitely-generated projectives, and which also have bounded homology.
So one has the inclusions $\Dc(R) \subseteq \Kpb{R} \subseteq \DR$.
Note that $\Kpb{R}$ is the image in $\DR$ of the Waldhausen
subcategory $\U_G(R) \subseteq \ch_R$.  It is an easy exercise to
check that when $R$ is right-Noetherian one has $\Kpb{R}=\D_b(\modh
R)$, where the latter denotes the full subcategory of $\DR$ consisting
of the bounded complexes of finitely-generated modules.

Theorem~\ref{thm-Gequiv} follows immediately from the following more
comprehensive statement:

\begin{theorem}
Let $R$ and $S$ be right-Noetherian.
\label{thm-Gequiv2}
\begin{enumerate}[(a)]
\item If $R$ and $S$ are derived equivalent, then $G(R)\he G(S)$.
\item $R$ and $S$ are
derived equivalent if and only if $\Kpb{R}$ and
$\Kpb{S}$ are equivalent as triangulated categories.
\item If $\D_b(\modh R)\he_{\Delta} \D_b( \modh S)$, then $G(R)\he
G(S)$ and $K(R)\he K(S)$.
\end{enumerate}
\end{theorem}

\begin{proof}
Part (b) is entirely due to Rickard \cite[8.1,8.2]{rickard-Morita}.
(Note that Rickard uses cochain complexes whereas we use chain
complexes, and writes $K^{-,b}(\proj R)$ for what we call
$\Kpb{R}$, etc.)  

For (a), suppose that $R$ and $S$ are derived equivalent.  Then
Theorem~\ref{thm-tilting one} says that there is a chain of Quillen
equivalences between $\ch_R$ and $\ch_S$.  On the homotopy categories,
this gives us a chain of triangulated equivalences between $\DR$ and
$\DS$.  Proposition~\ref{pr-DRsubcat} says that this triangulated
equivalence between $\DR$ and $\DS$ restricts to an equivalence
between $\Kpb{R}$ and $\Kpb{S}$.  So the complete Waldhausen
subcategory $\cU_G(R)$ is carried to $\cU_G(S)$ via the various
adjoint functors in the chain of Quillen equivalences.  One can now
use Corollary~\ref{co:QuilKiso2} to deduce that $K(\cU_G(R)) \he
K(\cU_G(S))$.  That is, $G(R)\he G(S)$.

For (c), recall that when $R$ is Noetherian $\D_b(\modh R)$ is just
another name for $\Kpb{R}$, and the same for $S$.  So if $\D_b(\modh
R)\he \D_b(\modh S)$ then by (b) $R$ and $S$ are derived equivalent; so
we can apply (a) and Theorem~\ref{thm-Kequiv}.
\end{proof}

%%%%%%%%%%%%%%%%%%%%%%%%%%%%%%%%%%%%%%%%%%%%%%%%%%%%%%%%%%%%%%%%%%%%%%%%

\section{Derived equivalence implies Quillen equivalence}
\label{sec-proof de iff qe}

In this section we prove the Tilting Theorem~\ref{thm-tilting one}.
The only difficult part of this theorem follows from a differential
graded analogue of the following result from~\cite[V.1]{gabriel}.  This can also
be viewed as another perspective on Morita theory.

\begin{theorem}{\bf (Gabriel)}\label{thm-gab}
Let $\cA$ be a co-complete abelian category with a small, projective,
strong generator $P$.  Then the functor
\[ \hom_{\cA}(P, -) \mc \cA \to \Mod\hom_{\cA}(P,P) \]
is the right adjoint of an equivalence of categories.
\end{theorem}

\noindent
There is also a version of this theorem for a set of small generators,
due to Freyd; see Section~\ref{sec-many}.

We begin by defining a chain complex of morphisms between any two
chain complexes.  For $M, N$ in $\ch_R$ define $\CHom(M,N)$ in $\chZ$
by
\[ \CHom(M, N)_n = \prod_k \hom_R(M_k, N_{n+k}).\] 
The differential for $\CHom(M,N)$ is given by $df_n = d_N f_n +
(-1)^{n+1} f_n d_M$.  This structure gives an enrichment of $\ch_R$
over $\chZ$.  So instead of an endomorphism {\it ring\/}, an object in
$\ch_R$ has a differential graded ring of endomorphisms.

\begin{definition}\label{def-dga}
{\em The tensor product of $X$ and $Y$ in $\chZ$ is defined by
$$(X \otimes Y)_n  = \bigoplus_k  X_k \otimes Y_{n-k}$$
where $d(x_p \otimes y_q) = dx_p \otimes y_q + (-1)^p x_p \otimes dy_q$.
A \dfn{differential graded algebra} is a chain complex $A$ in $\chZ$
with an associative and unital multiplication 
$\mu\mc A \otimes A \to A$~\cite[4.5.2]{weibel}.
A (right) \dfn{differential graded module} M over a 
differential graded algebra $A$ is a
chain complex $M$ with an associative and unital action $\alpha \mc
M \otimes A \to A$.  Denote the category of such modules by $\ModA$.  
}\end{definition}

For any $P$ in $\ch_R$ let $\EPo=\CHom(P,P)$.  Notice that $\EPo$ is a
differential graded ring with the product structure coming from
composition.  Also, for any $X\in \ch_R$ the complex $\CHom(P,X)$ is a
right differential graded $\EPo$-module with the action given by
precomposition.  So $\CHom (P,-)$ induces a functor from $\ch_R$ to
$\Mod\EPo$.  Its left adjoint is denoted $-\otimes_{\EPo} P$.  This
left adjoint can be defined as the coequalizer that the notation
suggests using the evaluation map $\CHom(P,P) \otimes P \to P$.

Our differential graded analogue of Gabriel's theorem produces a
Quillen equivalence of model categories instead of an equivalence of
categories.  So before stating it we need to establish the model
category structure on a category of differential graded modules.  The
following proposition is proved in~\cite[2.2.1, 3.1]{hinich} and
in~\cite[4.1.1]{ss1}.

\begin{proposition}\label{prop-modr mc}
Let $A$ be a DGA.  The category $\Mod A$ has a model category
structure where the weak equivalences are the maps inducing an
isomorphism in homology and the fibrations are the surjections.  The
cofibrations are then determined to be the maps with the
left-lifting-property with respect to the trivial fibrations.
\end{proposition}

We can now state the following differential graded version of
Gabriel's theorem.

\begin{theorem} 
\label{thm-gabriel one} Let 
$P$ in $\ch_R$ be a bounded complex of finitely generated projectives.
If $P$ is a (weak) generator for $\ch_R$, then there is a Quillen
equivalence
\[ \Mod\EPo  \ \varrow{1cm}\ {\ch_R} \]
in which the right-adjoint is the functor $\CHom(P,-)$.
\end{theorem}

Before proving this theorem we need the following lemma.

\begin{lemma}\label{lem-H}
Let $M, N\in \ch_R$.  Then
$H_*\CHom(M,N) \iso \DR(M,N)_*$ when $M$ is cofibrant. 
\end{lemma}

\begin{proof}
It is easy to see in general that $H_n\Hom_R(M,N)\iso H_n
\Hom_R(\Sigma^n M, N) \iso [\Sigma^n M,N]$ where $[-,-]$
denotes chain-homotopy-classes of maps.  When $A$ is cofibrant
one has that $\D_R(A,B)\iso [A,B]$ (since all objects are fibrant in
$\ch_R$), and so we can write
\[ H_n\Hom_R(M,N)\iso [\Sigma^n M,N]\iso \D_R(\Sigma^n M,N)
=\D_R(M,N)_n.
\]
\end{proof}

\begin{proof}[Proof of Theorem~\ref{thm-gabriel one}]
For any complex of projectives $P$, $\CHom(P,-)$ preserves surjections
(fibrations) and hence is exact.  We next show that $\CHom(P,-)$
preserves trivial fibrations; since $\CHom(P,-)$ is exact, we only
need to show that $H_*\CHom(P,K)=0$ when $H_*K = 0$ and apply this to
the kernel $K$ of the trivial fibration.  $P$ is cofibrant
by~\cite[2.3.6]{hovey-book} because $P$ is a bounded complex of
projectives.  Thus, by Lemma~\ref{lem-H}, if $K$ is acyclic then
$H_*\CHom(P,K) \iso \DR(P,K)_* \iso 0$.  Hence, the functor
$\CHom(P,-)$ preserves fibrations and trivial fibrations; see
also~\cite[4.2.13]{hovey-book}.  So its left adjoint is a Quillen map,
and therefore the adjoint pair induces total derived functors on the
level of homotopy categories~\cite[I.4]{Q}.  Denote these derived
functors by $\RCHom(P, -)$ and $ - \otimes^L_{\EPo} P$ respectively.

Since $\Ch_R$ and $\Mod\EPo$ are stable model categories, both total
derived functors preserve shifts and triangles in the homotopy
categories, i.e., they are exact functors of triangulated categories
by~\cite[I.4 Prop.\ 2]{Q}.  Since $- \otimes^L_{\EPo} P$ is a left
adjoint it commutes with coproducts.  To see that $\RCHom(P,-)$
commutes with coproducts it is enough to show that
$\DE(\EPo,\RCHom(P,-))$ commutes with coproducts since
$\EPo$ is a compact generator of $\Mod\EPo$.  By adjointness, this
functor is isomorphic to $\DR( \EPo \otimes^L_{\EPo} P, -)$ which in
turn is isomorphic to $\DR(P,-)$ since $\EPo$ is cofibrant.  
Since $P$ is
compact~\cite[6.4]{boek-nee} this functor commutes with coproducts.

Now consider the full subcategories of those $M$ in $\Ho(\Mod\EPo)$
and $X$ in $\DR$ respectively for which the unit of the adjunction
\[ \eta \ : \ M \ \varrow{1cm} \ \RCHom (P, M \otimes^L_{\EPo}  P) \]
or the counit of the adjunction
\[ \nu \ : \RCHom (P,X) \otimes^L_{\EPo} X \ \varrow{1cm} \ X \]
are isomorphisms. Since both derived functors are exact and preserve
coproducts, these are localizing subcategories.
The map $\eta$ is an isomorphism on the free module $\EPo$ and 
the map $\nu$ is an isomorphism on $P$.
Since the free module $\EPo$ generates the homotopy category
of $\EPo$-modules and $P$ generates $\Ch_R$, the derived
functors are inverse equivalences of the homotopy categories.
\end{proof}

Before completing the proof of the Tilting Theorem, here are two 
important statements.  

\begin{lemma}
\label{le-EMunique}
Suppose that $A$ is a DGA and $R$ is a ring (considered as a DGA
concentrated in degree zero).  Then $A$ and $R$ are quasi-isomorphic
if and only if $H_k(A)\iso H_k(R)$ for all $k$.  (That is, if and only
if $H_k(A) = 0$ for $k \neq 0$ and $H_0(A) \iso R$.)
\end{lemma}

\begin{proof}
Given $H_k(A) \iso H_k(R)$ for all $k$, then there are quasi-isomorphisms
of DGAs  $ A \xleftarrow{} A\langle 0\rangle \xrightarrow{} H_0(A) \iso R$.  
Here $A\langle 0\rangle $ is the (-1)-connected cover of $A$ with 
$A\langle 0\rangle_k = 0$ for $k < 0$, $A\langle 0\rangle_k = A_k$
for $k > 0$ and $A\langle 0\rangle_0 = Z_0 A$ the zero cycles.
\end{proof}

\begin{proposition}\label{qiso}
Any quasi-isomorphism $A \to B$ of differential graded algebras
induces a Quillen equivalence $\ModA \ra \Mod B$.
\end{proposition}

\begin{proof}
Any map $f \mc A \to B$ induces a Quillen adjoint pair between $\ModA$
and $\Mod B$, just as in Example~\ref{ex:Quilpair}.  The right adjoint
is given by restriction of scalars and the left adjoint is $-
\otimes_A B$.  ~\cite[4.3]{ss1} shows that this adjoint pair is a
Quillen equivalence.
\end{proof}

\begin{proof}[Completion of the proof of Theorem \ref{thm-tilting
one}] 
We must show that $(4)\Rightarrow (1)$, so suppose that $\ch_R$ has a
tilting complex $T$.  Then $T$ satisfies the hypotheses of Theorem
\ref{thm-gabriel one}, hence $\ch_R$ is Quillen equivalent to the
category of modules over the differential graded algebra $\EndR(T)$.
Since $T$ is a bounded complex of projectives, it is cofibrant
by~\cite[2.3.6]{hovey-book}; hence from Lemma~\ref{lem-H} we have
$H_*\EndR(T) \iso \DR(T,T)_*\iso S$ concentrated in dimension zero.
By Lemma~\ref{le-EMunique} this implies that $\EndR(T)$ is
quasi-isomorphic to $S$.  Thus the categories of $\EndR(T)$-modules
and right differential graded $S$-modules ($\ch_S$) are Quillen
equivalent by Proposition~\ref{qiso}:
\[ \ch_R \quad \simeq_Q \quad \Mod\End(T) \quad \simeq_Q \quad
\ch_S.   
\]
\end{proof}

%%%%%%%%%%%%%%%%%%%%%%%%%%%%%%%%%%%%%

\begin{remark}
\label{rem-marco}
We have now shown that when $R$ and $S$ are rings, their model
categories of dg-modules are Quillen equivalent if and only if the
associated homotopy categories are triangulated equivalent.  This is
false if $R$ and $S$ are allowed to be DGAs rather than rings,
essentially because the analog of Lemma~\ref{le-EMunique} fails: the
quasi-isomorphism type of an arbitrary DGA is not determined by its
homology (not even if you include all its Massey products,
see~\cite[A.3]{QT}).

In~\cite{ds-marco} we give an explicit example of two DGAs which are
derived equivalent, but where the model categories of dg-modules are
not Quillen equivalent.  The example is based on~\cite{marco} which
considers model categories underlying the stable category of modules
over the Frobenius rings $R = \bZ/p^2$ and $R' = \bZ/p
[\epsilon]/\epsilon^2$.  The homotopy categories are triangulated
equivalent but the corresponding $K$-theory groups are non-isomorphic
at $K_4$.  So by Corollary~\ref{co:stable QuilKiso} these model
categories cannot be Quillen equivalent.  In \cite{ds-marco} we give a
simpler proof of this by studying certain endomorphism DGAs, where we
can detect the difference in the second Postnikov sections instead of
in $K_4$.
\end{remark}

%%%%%%%%%%%%%%%%%%%%%%%%%%%%%%%%%%%%%%%%%%%%%%%%%%%%%%%%%%%%%%%%%%%%%%%

\section{Many generators version of proofs}\label{sec-many}

In this section we generalize the work in Section~\ref{sec-proof de iff qe}
to the case where we have a set of generators instead of just one.
Here the analogue for abelian categories is in~\cite[5.3H]{freyd}.
For derived equivalences, Keller~\cite[8.2]{dg-keller} gave the
corresponding extension of  Rickard's work~\cite[6.4]{rickard-Morita}. 
As always, our purpose is just to upgrade the derived equivalences to
Quillen equivalences.

\smallskip

As in Section~\ref{sec-proof de iff qe}, before moving to a
differential graded setting we first recall the classical setting.
Define a {\em ring with many objects} to be a small $\ab$-category (a
category enriched over abelian groups); this terminology makes sense
because an $\ab$-category with one object corresponds to a ring, with
composition corresponding to the ring multiplication.  Given a ring
with many objects $\uR$, a (right) $\uR$-module $M$ is a contravariant
additive functor from $\uR$ to $\ab$.  This means that for any two
objects $P, P'$ in $\uR$ there are maps $M(P') \otimes \uR(P,P') \to
M(P)$.  The category of right $\uR$-modules is an abelian category.

If $\cA$ is an abelian category and $\cP$ is a set of objects, say
that $\cP$ is a \dfn{set of strong generators} if $X=0$ whenever
$\hom_\cA(P,X)=0$ for every $P$ in $\cP$.  Define $\mEnd_{\cA}(\cP)$ to
be the full subcategory of $\cA$ (enriched over $\ab$) with object set
$\cP$.  The following theorem from~\cite[5.3H]{freyd} classifies
abelian categories with a set of strong generators:

\begin{theorem}{\bf (Freyd)}\label{thm-freyd}
Let $\cA$ be a co-complete abelian category with a set of small, projective,
strong generators $\cP$.  Then the functor
\[ 
\hom_{\cA}(P, -) \mc \cA \to \Mod\mEnd_{\cA}(\cP) \]
is the right adjoint of an equivalence of categories.
\end{theorem}

In order to generalize this result to a more homotopical setting, we
need to replace $\ab$-categories with $\ch$-categories (categories
enriched over $\ch = \ch_{\bZ}$.)  Since a $\ch$-category with one
object is a differential graded algebra, one may think of a small
$\ch$-category as a \dfn{DGA with many objects}.  Given a small
$\ch$-category $\cR$, a (right) $\cR$-module $\cM$ is a contravariant
$\ch$-functor from $\cR$ to $\ch$.  This means that for any two
objects $P, P'$ of $\cR$ there is a structure map of chain complexes
$\cM(P') \otimes \cR(P, P') \to \cM(P)$.  See~\cite[1.2]{kelly}
or~\cite[6.2]{bor} for more details.

Notice that $\ch_R$ and $\ch_{\uR}$ are both $\ch$-categories, where
$R$ is a ring and $\uR$ is a ring with many objects.  The enrichment
of $\ch_R$ over $\ch$ was discussed in the previous section.  Since
any two $\uR$-modules have an abelian group of morphisms
$\hom_{\uR}(M,N)$, the enrichment for $\ch_{\uR}$ follows similarly.

\begin{defn} \label{def-EP}
{\em Let $\cP$ be a set of objects in a $\ch$-category $\C$.
We denote by $\EP$ the full subcategory of $\C$ (enriched over $\ch$) 
with objects $\cP$, i.e., $\EP(P,P')=\HomC (P,P')$.   We let
\[ \HomC (\cP,-) \, : \,  \C \ \varrow{1cm} \ \Mod\EP \]
denote the  functor given by $\HomC (\cP,Y)(P) = \HomC (P,Y)$.}
\end{defn}

Note that if $\cP= \{ P \}$ has a single element, then $\EP$ is
determined by the single differential graded ring $\EndC (P) = \HomC
(P,P)$.  

In \cite[6.1]{emmc} it is established that there is a (projective) model
structure on the category $\Mod\EP$ of $\EP$-modules: the weak
equivalences are the maps which
induce quasi-isomorphisms at each object and the fibrations are the
epimorphisms (at each object).  

Now we can state the differential graded analogue of Freyd's theorem;
the difference is that here we have weak generators and a Quillen
equivalence instead of strong generators and a categorical
equivalence.  A set of objects $\cP$ in a stable model category $\C$
is a \dfn{set of (weak) generators} if the only localizing subcategory
of $\Ho(\C)$ which contains $\cP$ is $\Ho(\C)$.  As mentioned above
Theorem~\ref{thm-tilting one}, when the elements of $\cP$ are compact
then they generate $\Ho(\cC)$ if and only if they can detect when an
object is trivial; see~\cite[2.2.1]{ss2}.  Note that a (possibly
infinite) coproduct of a set of generators is still a generator, but
is not necessarily compact.

\begin{theorem} \label{thm-freyd1} Let $\uR$
be a ring with many objects and  $\cP$ a set in $\ch_{\uR}$
of bounded complexes of finitely generated projectives. 
If $\cP$ is a set of (weak) generators for $\ch_{\uR}$ then 
there is a Quillen equivalence
\[ \Mod\EP \ \varrow{1cm} \ {\ch_{\uR}} \]
in which the right adjoint is the functor $\Hom_{\uR}(\cP,-)$.
\end{theorem}

Note that for every object $r\in \uR$ there is a corresponding {\em
`free module'} $F_r^{\uR}$ given by $F_r^{\uR}(s)=\uR(s,r)$.
A projective $\uR$-module is {\it finitely-generated\/} if it is a direct
summand of a module $\oplus_i F_{r_i}^{\uR}$, where the sum is finite.
And as usual, we denote the homotopy category of $\ch_{\uR}$ by $\D_{\uR}$.
We need the following lemma:

\begin{lemma}
\label{le-keller}
The compact objects in $\D_{\uR}$ are those complexes which are
quasi-isomorphic to a bounded complex of finitely-generated
projective $\uR$-modules.
\end{lemma}

\begin{proof}
This follows from~\cite[5.3]{dg-keller}.
\end{proof}

\begin{proof}[Proof of Theorem~\ref{thm-freyd1}]
Just as in the proof of Theorem~\ref{thm-gabriel one}, one can check
that $\Hom_{\uR}(P, -)$ takes fibrations and trivial fibrations in
$\ch_{\uR}$ to fibrations and trivial fibrations in $\ch$ for any
bounded complex of projectives $P$.  So the functor $\Hom_{\uR}(\cP,
-)$ preserves fibrations and trivial fibrations.  Thus, together with
its left adjoint $-\otimes_{\EP}\cP$, it forms a Quillen pair.

We proceed as in the proof of Theorem~\ref{thm-gabriel one}.  The
induced total derived functors are again exact functors of
triangulated categories which commute with coproducts.  Here we use
the fact that each $P$ is compact to show the right adjoint commutes
with coproducts.  The full subcategories for which the unit of the
adjunction $\eta$ or the counit of the adjunction $\nu$ are
isomorphisms are localizing subcategories.  

Note that for each object $P$ in $\cP$ there is a free $\EP$-module
$F_P^{\EP}$ defined by $F_P^{\EP}(P')=\EP(P',P)$, and these generate
the homotopy category of $\EP$-modules.  For every $P\in\cP$ the
$\EP$-module $\Hom_{\uR} (\cP,P)$ is isomorphic to the free module
$F_P^{\EP}$ by inspection, and $F_P^{\EP} \otimes_{\EP} \cP$ is
isomorphic to $P$ since they represent the same functor on
$\ch_{\uR}$.  Thus, $\eta$ is an isomorphism on every free module and
$\nu$ is an isomorphism on every object of $\cP$.  Since the free
modules $F_P^{\EP}$ generate the homotopy category of $\EP$-modules
and the objects of $\cP$ generate $\ch_{\uR}$, the localizing
subcategories where $\eta$ and $\nu$ are isomorphisms are the whole
homotopy categories.  This implies that the adjoint pair is a Quillen
equivalence.
\end{proof}

Finally, we can write down a many-objects version of
Theorem~\ref{thm-tilting one}.  If $\cP$ is a set of (weak) generators
with each element $P$ a bounded complex of finitely generated
projectives and $H_*\EP$ is concentrated in degree zero, then we call
$\cP$ a set of {\em tiltors}.  The following theorem is a
generalization of Keller's work~\cite[8.2]{dg-keller}:

\begin{theorem} \label{thm-tilting} {\bf (Many-objects tilting theorem)}
Theorem~\ref{thm-tilting one} holds when the rings $R$ and $S$ are
replaced by rings-with-many-objects $\uR$ and $\uS$.  The tilting
complex is replaced by a set of tiltors $\bT$ with $H_*\cE(\bT) \iso
\uS$.
\end{theorem}

The proof is given below, but first we state some easy consequences:

\begin{corollary}\label{cor-ringoids}
Two rings-with-many-objects $\uR$ and $\uS$ are derived equivalent if
and only if their associated model categories of chain complexes
$\ch_{\uR}$ and $\ch_{\uS}$ are Quillen equivalent.
\end{corollary}

Using Theorem~\ref{thm-freyd} we get the following corollary as well.
Given an abelian category $\cA$ satisfying the hypotheses of
Theorem~\ref{thm-freyd}, choose a set of small, projective, strong
generators $\cP$.  Let $\uA = \mEnd_{\cA}(\cP)$ be the associated
ring-with-many-objects.  Freyd's theorem says that $\cA$ is equivalent
to $\Mod\uA$, and so $\ch_{\cA}$ is equivalent to $\ch_{\uA}$.  In
particular, one gets a projective model structure on $\ch_{\cA}$ by
lifting the one on $\ch_{\uA}$ across the equivalence;
see~\cite[6.1]{emmc}.  The next result is now an immediate consequence
of Corollary~\ref{cor-ringoids}.

\begin{corollary}\label{cor-abcat}
Let $\cA$ and $\cB$ be co-complete abelian categories with sets of small,
projective, strong generators.  Then $\cA$ and $\cB$ are derived
equivalent if and only if their associated model categories of chain
complexes $\ch_{\cA}$ and $\ch_{\cB}$ are Quillen equivalent.
\end{corollary}

\begin{remark}\label{caution}
Warning: Let $\cM$ and $\cN$ be two stable model categories whose
underlying categories are abelian, with sets of small, strong,
projective generators.  The above corollary does {\it not\/} say that
$\cM$ and $\cN$ are Quillen equivalent if and only if $\Ho(\cM)$ and
$\Ho(\cN)$ are triangulated equivalent.  This statement is false; see
~\cite{marco},~\cite{ds-marco}.  Note in particular that it does not
apply to the model category $\Mod\cR$ where $\cR$ is a DGA; the
problem is that $\Ho(\Mod\cR)$ is not the same as
$\Ho(\ch_{\Mod\cR})$.
\end{remark}

\begin{proof}
[Proof of Theorem~\ref{thm-abresult}] Part (a) is the above corollary.
Part (b) is immediate from (a) and Corollary~\ref{co:stable QuilKiso}.
\end{proof}

\begin{proof}
[Proof of Theorem \ref{thm-tilting}]
The proof that condition (1) implies condition (2)
and condition (2) implies condition (3) follows just as in 
Theorem~\ref{thm-tilting one}. 

Now assume condition (3) and fix a triangulated equivalence between
$K_b(\mbox{proj-}\uR)$ and $K_b(\mbox{proj-}\uS)$.  For $s$ any object
in $\uS$, consider the module $F_s^{\uS}$ as a complex concentrated in
dimension zero; let $T_s$ be its image in $K_b(\mbox{proj-}\uR)$.
Since the objects in $\{F_s\}_{s \in \uS}$ generate
$K_b(\mbox{proj-}\uS)$, the objects in $\bT =\{T_s\}_{s\in\uS}$
generate $K_b(\mbox{proj-}\uR)$.  But $K_b(\mbox{proj-}\uR)$
generates $\D_R$, so $\bT$ generates $\D_R$ as well.  By
Lemma~\ref{le-keller} the objects of $\bT$ are compact in $\D_R$.
Finally, we
also have $H_*\cE(\bT) \iso H_*\cE(\{ F_s^{\uS}\}) \iso\uS$.  So $\bT$ is a
set of tiltors.

If we are given a set of tiltors $\mathbb T$ for $\ch_{\uR}$, then by
Theorem~\ref{thm-freyd1}  $\ch_{\uR}$ is Quillen equivalent to the
category of modules over the endomorphism category $\cE(\mathbb T)$.
Since $H_*\cE(\mathbb T)\iso \uS$ is concentrated in dimension
zero, $\cE(\mathbb T)$ is quasi-isomorphic to $\uS$ by an extension of
Lemma~\ref{le-EMunique}.  Thus the categories of differential graded
$\cE(\mathbb T)$-modules and differential graded $\uS$-modules are
Quillen equivalent by~\cite[6.1]{emmc} which generalizes
Proposition~\ref{qiso}:
\[ \ch_{\uR} \quad \simeq_Q \quad \Mod\cE(\mathbb T) \quad \simeq_Q \quad
\ch_{\uS}   \]
\end{proof}

%%%%%%%%%%%%%%%%%%%%%%%%%%%%%%%%%%%%%%%%%%%%%%%%%%%%%%%%%%%%%%%%

\appendix\label{app}

\section{Proof of Proposition~\ref{pr:QuilKiso}}

Recall that $\cM$ and $\cN$ are pointed model categories, $(L,R)\colon
\cM \ra \cN$ is a Quillen equivalence, $\cU$ is a complete Waldhausen
subcategory of $\cM$, and $\cV=\overline{(L\cU)}$.  (Note that $L$,
being a left adjoint, must preserve the initial object).  We must show
that $\cV$ is a complete Waldhausen subcategory of $\cM$ and that the
induced map $L\colon K(\cU) \ra K(\cV)$ is a weak equivalence.

For the remainder of this section, let $F$ be a fibrant-replacement
functor in $\cN$ and let $Q$ be a cofibrant-replacement functor in
$\cM$.  Note that the functor $QRF\colon \cN \ra \cM$ takes $\cV$ into
$\cU$: for if $X\in \cV$ then $X\he LA$ for some $A\in \cU$, and then
$QRFX\he QRFLA \he A$.  Since $\cU$ is complete and $A\in \cU$, it
follows that $QRFX\in \cU$ as well.  

\begin{lemma}
\label{le-klem1}
$\cV$ is a complete Waldhausen subcategory of $\cN$.
\end{lemma}

\begin{proof}
The only point which takes work is axiom (iii) for Waldhausen
categories.  So if $A\cofib B$ and $A\ra X$ are maps in $\cN$ where
$A$, $B$, $X \in \cV$, we must show that the pushout
$B\amalg_A X$ is also in $\cV$.

Consider the maps $QRFA \ra QRFB$ and $QRFA \ra QRFX$.  All the
domains and codomains of these maps are in $\cU$.  Factor $QRFA\ra
QRFB$ as $QRFA \cofib Z \trfib QRFB$.  Then $Z\in \cU$ and so the
pushout $P=Z\amalg_{QRFA} QRFX$ is also in $\cU$, because $\cU$ is a
Waldhausen subcategory of $\cM$.  This pushout is weakly equivalent to
the homotopy pushout (see~\cite[10]{DS}) of $Z \la QRFA \ra QRFX$,
because $QRFA \ra Z$ is a cofibration and all the objects $Z$, $QRFA$,
and $QRFX$ are cofibrant; see~\cite[5.2.6]{hovey-book}.  Since $Z
\trfib QRFB$, $P$ is also weakly equivalent to the homotopy pushout of
the diagram $QRFB\la QRFA \ra QRFX$.

Finally, any left Quillen functor $L$ preserves homotopy pushouts, in
the sense that $LP$ is weakly equivalent to the homotopy pushout of
$LQRFB \la LQRFA \la LQRFX$.  The latter homotopy pushout is weakly
equivalent to the homotopy pushout of $B\la A \ra X$, which in turn is
just weakly equivalent to the pushout $B\amalg_A X$ (since $A\ra B$ is
a cofibration and all the objects $A$, $B$, $X$ are cofibrant).  So
$B\amalg_A X$ is weakly equivalent to $LP$, and is therefore in
$\overline{L\cU}$.
\end{proof}

Let $w\cU$ denote the subcategory consisting of all weak equivalences
in $\cU$, and write $N(w\cU)$ for the nerve of this category.  The
functor $L$ induces a map $w\cU \ra w\cV$.

\begin{lemma}
\label{le-klem2}
$NL\mc N(w\cU) \ra N(w\cV)$ is a weak equivalence of spaces.
\end{lemma}

\begin{proof}
First of all, the functor $Q\colon \cM \ra \cM$ maps $\cU$ into itself
(because $\cU$ is complete), and comes equipped with a natural
transformation $QX\ra X$.  This shows that the induced map $NQ\colon
N\cU \ra N\cU$ is homotopic to the identity~\cite[2.1]{segal}.
Similarly, $NF\colon N\cV \ra N\cV$ is homotopic to the identity.

The functor $QRF\colon \cN \ra \cM$ maps $\cV$ to $\cU$, as was
remarked prior to the previous lemma.  There are natural
transformations $LQRF \ra LRF \ra F$, and $Q\ra QRL \ra QRFL$.
It follows that the compositions $NL\circ N(QRF)$ and $N(QRF)\circ
NL$ are homotopic to the respective identity maps, and so are part of
a homotopy equivalence.    
\end{proof}

Let $\Delta_n$ denote the category consisting of $n$ composable arrows
$0 \ra 1 \ra \cdots \ra n$.  This may be given the structure of a {\it
Reedy category\/}~\cite[5.2.1]{hovey-book} in which all the maps
increase dimension.  The category of diagrams $\cM^{\Delta_n}$ has a
corresponding {\it Reedy model structure\/} \cite[5.2.5]{hovey-book}
in which a map $X_\wdo \ra Y_\wdo$ is a weak equivalence (respectively
fibration) if and only if each $X_n \ra Y_n$ is a weak equivalence
(respectively fibration).  A map is a cofibration if and only if all
the maps $X_n\amalg_{X_{n-1}} Y_{n-1} \ra Y_n$ are cofibrations.  In
particular, an object $X_\wdo$ is cofibrant if and only if the maps
$X_{n-1} \ra X_n$ are all cofibrations; by a simple recursion, this
implies that all the $X_i$'s are cofibrant as well.

Let $\cU_n$ denote the full subcategory of $\cM^{\Delta_n}$ consisting
of cofibrant diagrams whose objects all belong to $\cU$.  It is easy
to see that $\cU_n$ is a complete Waldhausen subcategory of
$\cM^{\Delta_n}$.  The functors $(L,R)$ prolong to functors
$(L,R)\colon \cM^{\Delta_n} \ra \cN^{\Delta_n}$, and this is still a
Quillen equivalence.  We need the following:

\begin{lemma}
Any diagram in $\cV_n$ is weakly equivalent to one in
$L(\cU_n)$---i.e., $\cV_n=\overline{L\cU_n}$.
\end{lemma}

\begin{proof}
Let $X_\wdo=[X_0\ra \cdots \ra X_n]$ be an object in $\cV_n$.  Then
each $X_i$ is in $\cV$, and so $QRFX_i$ lies in $\cU$ (as was shown
above Lemma~\ref{le-klem1}).  Consider the object $QRFX_\wdo=[QRFX_0
\ra \cdots \ra QRFX_n]$.  This need not be cofibrant in
$\cM^{\Delta_n}$, but we can still take its cofibrant
replacement---call this new object $C_\wdo$.  Each $C_i$ is a
cofibrant object weakly equivalent to $QRFX_i$, and is therefore in
$\cU$; so $C_\wdo$ is in $\cU_n$.  We have a sequence of maps $LC_\wdo
\ra LQRFX_\wdo \ra FX_\wdo \la X_\wdo$, all of which are objectwise
weak equivalences, and so $X_\wdo$ is weakly equivalent to an object
in $L\cU_n$.
\end{proof}

The category $w(\cU_n)$ is exactly the category $wF_{n}(\cU)$ defined
in Section~\ref{se:Kthy}.  So there is a `forgetful' functor
$wS_n(\cU) \ra w(\cU_{n-1})$: in the notation of
\cite[1.3]{waldhausen} it sends an object $\{A_{ij}\}$ to the sequence
$A_{01} \cofib A_{02} \cofib \cdots \cofib A_{0n}$.  This functor is
easily seen to be an equivalence of categories (see \cite[bottom of
p. 328]{waldhausen}).

\begin{proof}[Proof of Proposition~\ref{pr:QuilKiso}]
Recall that $K(\cU)$ is defined as the geometric realization of a
simplicial space $[n] \mapsto N(wS_n(\cU))$.  It is therefore 
enough to show that $L$ induces weak equivalences $N(wS_n(\cU)) \ra
N(wS_n(\cV))$ at each level.  There is a commutative diagram
\[ \xymatrix{ wS_n(\cU) \ar[d] \ar[r]  & wS_n (\cV) \ar[d] \\
              w\cU_{n-1} \ar[r] & w\cV_{n-1}
}
\]
and the vertical maps are equivalences of categories.  So it suffices
to show that the maps $N(w\cU_n) \ra N(w\cV_n)$ are weak equivalences.
But this follows from Lemma~\ref{le-klem2} applied to the complete
Waldhausen subcategories $\cU_n$ and $\cV_n=\overline{L\cU_n}$ of
$\cM^{\Delta^n}$.
\end{proof}

%%%%%%%%%%%%%%%%%%%%%%%%%%%%%%%%%%%%%%%%%%%%%%%%%%%%%%%%%%%%%%%%%%%

\end{document}